\documentclass{article}

\usepackage{proof}
\usepackage{latexsym}
\usepackage{amssymb}

\newtheorem{theorem}{Theorem}
\newtheorem{definition}[theorem]{Definition}
\newtheorem{proposition}[theorem]{Proposition}
\newtheorem{lemma}[theorem]{Lemma}

\begin{document}

\title{nested PLS}




\author{Toshiyasu Arai
\\
Graduate School of Science,
Chiba University
\\
1-33, Yayoi-cho, Inage-ku,
Chiba, 263-8522, JAPAN
}

\date{}

\maketitle

\begin{abstract}
In this note we will introduce a class of search problems, called {\it nested Polynomial Local Search\/}(nPLS)
problems, and show that definable NP search problems, i.e., $\Sigma^{b}_{1}$-definable functions 
 in $T^{2}_{2}$ are characterized in terms of the nested PLS.


\end{abstract}

\section{Introduction}\label{sect:intro}

TFNP defined in \cite{TFNP} is the class of searching problems in which a witness always exists
with polynomial time verification.
Papadimitriou et. al\cite{PLS}, \cite{Papa} introduced several subclasses of TFNP.
The idea is to group together problems in TFNP in terms of "proof styles",
where by "proof styles" a graph-theoretic lemma to prove the {\it totality\/}
of the problems is meant, e.g., a parity argument.

Consider the language of the bounded arithmetic with function symbols for polynomial time
computable functions, i.e., the language of the Cook's PV in the predicate logic.
Now a searching problem in TFNP can be described by a $\Sigma^{b}_{0}$-formula
$R(x,y)$ in the language such that $\forall x\exists y<2^{p(|x|)}\, R(x,y)$ is true for a polynomial $p$.
Hence by identifying "proof styles" with formal provability in a bounded arithmetic $T$,
a subclass $\mbox{TFNP}(T)$ of TFNP is associated with each $T$.
Then a $\Sigma^{b}_{1}$-definable (multivalued) function in $T$ computes
witnesses of a problem in $\mbox{TFNP}(T)$ and vice versa.

In \cite{BussKrajicek} Buss and Kraj\'\i\v cek proved that 
$\mbox{TFNP}(T^{1}_{2})=\mbox{PLS}$, {\it Polynomial Local Search problems\/}.
Let us see some details of the class PLS
since our new class of searching problems is an extension of it.

A problem ${\cal P}$ of PLS consists in the following:
\begin{enumerate}
\item

A polynomial time predicate $F(x,s)$  defines the set
\[
F(x)=\{s: F(x,s)\}
\]
of {\it feasible points\/}(search space) 
such that 
$s\in F(x)\to |s|\leq d(|x|)$
with a polynomial bound $d$.

\item
A polynomial time computable {\it initial point\/} function $i(x)\in F(x)$.
\item

A polynomial time computable {\it neighborhood\/} function $N(x,s)$ such that
\[
s\in F(x) \to N(x,s)\in F(x)
.\]
\item

A polynomial time computable {\it cost\/} function $c(x,s)$ such that
\[
N(x,s)=s \lor c(x,N(x,s))<c(x,s)
.\]
\end{enumerate}

Then an $s$ is said to be a solution or a goal of the instance $x$ iff
$s\in F(x) \land N(x,s)=s$.

Assuming that the problem is a PLS in $S^{1}_{2}$, i.e.,
that these conditions $s\in F(x)\to |s|\leq d(|x|)$, $i(x)\in F(x)$, $s\in F(x) \to N(x,s)\in F(x)$ and 
$N(x,s)=s \lor c(x,N(x,s))<c(x,s)$ are provable in $S^{1}_{2}$,
it is plain to see that $T^{1}_{2}$ proves its totality
\[
T^{1}_{2}\vdash \forall x \exists s<2^{d(|x|)}[s\in F(x) \land N(x,s)=s]
\]
since $T^{1}_{2}$ proves the existence of an $s\in F(x)$ with a (globally) minimal cost.

Note that in place of a neighborhood {\it function\/} $N(x,s)=t$,
a neighborhood {\it predicate\/} $N(x,s,t)$ yields a searching problem in PLS
if $N(x,s,t)$ is polynomial time computable and
\begin{eqnarray}
&&
\mbox{the cardinality of the set } \{t\in F(x): N(x,s,t)\}
\nonumber \\
&&
 \mbox{is bounded by a polynomial }
p(|x|) \mbox{ for any } s\in F(x).
\label{eq:cardinalitypolybnd}
\end{eqnarray}
For the predicate and a cost function $c$, 
an $s\in F(x)$ is said to be a solution of the instance $x$ iff
$s$ is a local minimum, i.e.,
$\forall t\in N(x,s)[c(x,s)\leq c(x,t)]$, which is polynomial time verifiable.

Let $N^{\prime}$ be a polynomial time computable predicate such that
\[
N^{\prime}(x,s,t) \Leftrightarrow
[s\neq t\land N(x,s,t)]\lor [s=t \land s \mbox{ is a solution}]
.\]
Then $s$ is a solution iff $N^{\prime}(x,s,s)$.

The graph-theoretic lemma invoked in a PLS problem is:
\begin{quote}
Every finite directed acyclic graph(DAG) has a sink.
\end{quote}
A DAG $(V,E)$ is associated with an instance $x$ of a PLS problem $(F,i,N,c)$ with a neighborhood predicate $N$ such that
the set of nodes $V=F(x)$, and $s$ is adjacent to $t$, $E(s,t)$ iff
$N(x,s,t)\land c(x,s)>c(x,t)$.
Then a sink of this DAG is a local minimum for the instance $x$.

\smallskip
Next let us consider $\mbox{TFNP}(T^{2}_{2})$.
Some characterizations of $\Sigma^{b}_{1}$-definable functions in $T^{2}_{2}$ 
have already known.

Kraj\'\i\v cek,  Skelly and Thapen\cite{KraSkeTha} characterized the class in terms of colored PLS problems.
Skelly and Thapen\cite{SkeTha} subsequently gave a characterization of the $\Sigma^{b}_{1}$-definable functions in $T^{k}_{2}$ for any $k\geq 2$, based on a combinatorial principle for $k$-turn games.

On the other side, Beckmann and Buss\cite{BeckmannBuss}, \cite{BeckmannBussnotation} gives a characterization of the $\Sigma^{b}_{i}$-definable functions of $T^{k+1}_{2}$ for all $0<i\leq k+1$, using a relativized notion of polynomial local search problems, $\Pi^{p}_{k}$-PLS problems with $\Pi^{p}_{i-1}$-goals.
As in PLS, 
a problem ${\cal P}$ of this class consists in a predicate $F(x,s)$ for feasible points,
an initial point function $i(x)$, a neighborhood function $N(x,s)$ and a cost function $c(x,s)$.
As contrasted with PLS, $F(x,s)$ is here a $\Pi^{p}_{k}$-predicate, and
a {\it goal\/} predicate such that
\[
G(x,s)\leftarrow s\in F(x) \land N(x,s)=s
\]
has to be a $\Pi^{p}_{i-1}$-predicate.

Each instance ${\cal P}(x)$ has at least one solution $s\in F(x)$ such that $N(x,s)=s$.
Thus ${\cal P}$ defines a multifunction associating a solution $y$ of an instance ${\cal P}(x)$ to
the input $x$. It is written $y={\cal P}(x)$.

When $k=0, i=1$, the class of $\Pi^{p}_{0}$-PLS problems with $\Pi^{p}_{0}$-goals
is equivalent to the PLS problems.

In \cite{BeckmannBuss} the price to pay in characterizing the $\Sigma^{b}_{1}$-definable functions in 
higher fragments $T^{k}_{2}$ was to raise the complexities of the sets $F(x)$ of feasible points, 
but the predicate $s\in G(x)$ to be goals is still in P.

\smallskip
In this note we will introduce a class of search problems, called {\it nested Polynomial Local Search\/}(nPLS) problems.
Let us illustrate an instance of an nPLS in terms of graphs.
Let $(V,E)$ be a non-empty finite directed graph with a function $c:V\to\mathbb{N}$ such that
\begin{equation}\label{eq:costnPLS}
E(s,t) \land s\neq t\Rightarrow c(s)>c(t)
\end{equation}
Then the graph has no cycle except trivial ones $E(s,s)$.

Now assume that
\begin{equation}\label{eq:totalitynPLS}
\forall s\in V\exists t\in V[E(s,t)]
\end{equation}
Then it is straightforward to see that the graph has a trivial cycle.
Note here that the size of the set of neighborhoods $\{t: E(s,t)\}$ may be
exponentially large, cf. (\ref{eq:cardinalitypolybnd}).

However in searching a trivial cycle, another searching problems (\ref{eq:totalitynPLS}) are hidden to search a neighborhood $t$ to $s$.
Let us imagine that a witness of a searching problem $(V_{s},E_{s})$ 
associated with $s$ yields a neighborhood $t$, 
where $(V_{s},E_{s})$  is again a finite directed graph with a function $c_{s}$
enjoying (\ref{eq:costnPLS}) and (\ref{eq:totalitynPLS}), and a trivial cycle in $E_{s}$ is a witness sought.
Iterating this process, we get an nested PLS.
To avoid an infinite descent, a rank $rk(s)$ is associated with each
problem $s$ such that $((V_{s})_{t},(E_{s})_{t})$
is in lower rank than $(V_{s},E_{s})$, where
 $((V_{s})_{t},(E_{s})_{t})$ denotes the problem generated by $t\in V_{s}$,
 and to search a neighborhood to $t$ in $E_{s}$.
If $rk(s)=0$, then the problem $(V_{s},E_{s})$ is a PLS with (\ref{eq:cardinalitypolybnd}).

\smallskip
In the next section \ref{sect:nPLS} we will define the nPLS formally, and show that
$T^{2}_{2}$ proves its totality, Lemma \ref{lem:total}.

In the section \ref{sect:search} it is shown that $\mbox{TFNP}(T^{2}_{2})$ is an nPLS, Theorem \ref{th:nPLST22}.
Before proving the result, 
we will reprove the Buss-Kraj\'\i\v cek's result
$\mbox{TFNP}(T^{1}_{2})\subset\mbox{PLS}$ in \cite{BussKrajicek} 
to illuminate our proof method for searching a witness of a $\Sigma^{b}_{1}$-formula
in a derivation of bounded arithmetic.

\section{nested PLS}\label{sect:nPLS}

In this section we will introduce a class of search problems, called {\it nested Polynomial Local Search\/}(nPLS) problems.
In an nPLS the feasible points $s\in F(x)$ and the cost function are polynomial time computable,
while the neighborhood function is not a polynomial time computable function, but it's a multifunction
defined by an nPLS problem.

Let us call a PLS $0${\it -nPLS\/} (nPLS of rank $0$).
Assuming the class of $a$-nPLS problems has been defined for an integer $a$.
Then an $(a+1)$-nPLS consists in a polynomial time computable predicate $F(x,s)$,
polynomial time computable functions $i(x)$, $c(x,s)$ and 
 $a$-nPLS problems ${\cal P}_{x}(s)$ for $s\in F(x)$.
 
 $F(x,s)$ and $i(x)$ are as in PLS.
 Now each $s\in F(x)$ generates an $a$-nPLS ${\cal P}_{x}(s)$, which consists in
 a neighborhood {\it predicate\/} ${\cal N}_{x,s}(y,z)$, a result extracting function $U_{x,s}(y)$ such that
 
\[
{\cal N}_{x,s}(y,y) \to {\cal N}_{x}(s, U_{x,s}(y))
\]
where ${\cal N}_{x}(s,t)$ is a neighborhood predicate for the whole $(a+1)$-nPLS problem.
Namely, given a feasible point $s\in F(x)$, search a solution $y$ of the $a$-nPLS ${\cal P}_{x}(s)$, and then
extract result $U_{x,s}(y)$ from the solution $y$ to get a neighborhood to $s$.
 
Note that any $a$-nPLS is a $b$-nPLS for $b<a$
since for each polynomial time computable function $f$ there is a PLS ${\cal P}$
such that $f(x)={\cal P}(x)$, cf. \cite{BeckmannBuss}.
Moreover observe that for each fixed integer $a$,
a uniform $a$-nPLS problem is a PLS.

By unifying the iterated process of defining $a$-nPLS, we get an nPLS.

\begin{definition}\label{df:nPLS}
{\rm An} nPLS {\rm problem} ${\cal P}=\langle d; {\cal S},{\cal T},{\cal N}; N, i,t,c,S,U,rk\rangle$
{\rm where}
\begin{enumerate}
\item $d(x)$ {\rm is a polynomial.}

\item
${\cal S}(x,s)$ {\rm and} ${\cal T}(x,s,t)$ {\rm are  polynomial computable predicates for} 
sources {\rm and} targets{\rm , resp. with a polynomial bound}
$d$.

${\cal S}(x)=\{s: {\cal S}(x,s)\}$ {\rm and}
${\cal T}(x)=\sum\{{\cal T}(x,s): s\in{\cal S}(x)\}$ {\rm with} ${\cal T}(x,s)=\{t: {\cal T}(x,s,t)\}$ 

${\cal N}(x,s,y,z)$ {\rm is a polynomial time computable predicate for}
neighborhoods.

\item
$N(x,s,y)$, $i(x)$, $t(x,s)$, $c(x,s)$, $S(x,s,y)$, $U(x,s,y,z)$ {\rm and} $rk(x,s)$ {\rm are
polynomial time computable functions for} neighborhoods of rank zero, 
initial source, initial target,
cost, generated source, extracted result {\rm and} rank{\rm , resp.}

\end{enumerate}

{\rm These have to enjoy the following conditions.}
\begin{enumerate}

\item
\begin{equation}\label{eq:nPLS1}
s\in {\cal S}(x) \lor t\in{\cal T}(x,s) \to\max\{ |s|, |t|\}\leq d(|x|)
\end{equation}

\begin{equation}\label{eq:nPLS2.1}
y\in{\cal T}(x,s) \to 
S(x,s,y)\in{\cal S}(x)
\end{equation}

\begin{equation}\label{eq:nPLS2.4}
{\cal N}(x,s,y,z)  \to s\in{\cal S}(x)\land y,z\in{\cal T}(x,s)
\end{equation}

\item
\begin{equation}\label{eq:nPLS0}
rk(x,s)=0  \to 
[{\cal N}(x,s,y,z)\leftrightarrow N(x,s,y)=z\land y\in {\cal T}(x,s)]
\end{equation}

\item
\begin{equation}\label{eq:nPLS7}
rk(x,s)>0 \land y\in {\cal T}(x,s) \to {\cal N}(x,s,y,y) \lor rk(x,S(x,s,y))<rk(x,s)
\end{equation}

\item
\begin{equation}\label{eq:nPLS5}
rk(x,s)>0 \land
{\cal N}(x,S(x,s,y),z,z) 
\to 
{\cal N}(x,s,y,U(x,s,y,z))
\end{equation}

\item
\begin{equation}\label{eq:nPLS11}
i(x)\in {\cal S}(x)
\end{equation}

\begin{equation}\label{eq:nPLS12}
t(x,s)\in {\cal T}(x,s)
\end{equation}

\item
\begin{equation}\label{eq:nPLS13}
{\cal N}(x,s,y,z)\to y=z \lor c(x,y)>c(x,z)
\end{equation}

\end{enumerate}

{\rm If all of these conditions to be an nPLS are provable in} $S^{1}_{2}${\rm , then we say that the problem is an} nPLS in $S^{1}_{2}$.
\end{definition}

Let ${\cal P}$ be an nPLS. Then we write
\[
y={\cal P}(x) :\Leftrightarrow {\cal N}(x,i(x),y,y)
.\]
Let us depict a table showing structures of an nPLS.
\newpage
\begin{table}[h]
\caption{nPLS}
\label{tab:1}       
\begin{tabular}{|c|c|}
\hline\noalign{\smallskip}
sources  & targets
\\
\noalign{\smallskip}\hline\noalign{\smallskip}
$i(x)$ & 
\\
\noalign{\smallskip}\hline\noalign{\smallskip}
$s$ & $t(x,s)\mbox{\hspace{5mm}} y\to_{{\cal N}}U(x,s,y,z)\to_{{\cal N}}\cdots \to_{{\cal N}}u$
\\
\noalign{\smallskip}\hline\noalign{\smallskip}
$S(x,s,y)$ & solution $z\to_{{\cal N}}z$  
\\
\noalign{\smallskip}\hline\noalign{\smallskip}
\end{tabular}
\end{table}
\noindent
where each row $s$ denotes a search problem in the search space ${\cal T}(x,s)$, and
the higher row indicates search problems in higher ranks.
$y\to_{{\cal N}}z$ in the row $s$ designates ${\cal N}(x,s,y,z)$.

Given a $y\in {\cal T}(x,s)$, generate a lower rank source(problem) $S(x,s,y)$ and search a solution $z$ of the problem $S(x,s,y)$.
Then you will get a neighborhood $U(x,s,y,z)$ of $y$ in a lower cost.
Continue this search to find a solution $u$ of the problem $s$.
And if $s$ is not the highest rank problem, then extract a result from the solution $u$ to find a neighborhood of a higher rank target,
and so forth.

\begin{lemma}\label{lem:total}
Let ${\cal P}=\langle d; {\cal S},{\cal T},{\cal N};N,  i,t,c,S,U,rk\rangle$ be an nPLS in $S^{1}_{2}$.
Then
\[
T^{2}_{2}\vdash \forall x \exists y<2^{d(|x|)}[y={\cal P}(x)]
.\]
\end{lemma}
{\it Proof}.\hspace{2mm}
Argue in $T^{2}_{2}$.
We show by induction on $a$ that
\begin{equation}\label{eq:total1}
\forall s<2^{d(|x|)}\exists y<2^{d(|x|)} [s\in {\cal S}(x) \land rk(x,s)\leq a \to {\cal N}(x,s,y,y)]
\end{equation}

First consider the case $a=0$.
Then by (\ref{eq:nPLS2.4}),
(\ref{eq:nPLS0}) and (\ref{eq:nPLS13}), 
(\ref{eq:total1}) is equivalent to the totality of a PLS problem
\[
\forall s<2^{d(|x|)}\exists y<2^{d(|x|)} [s\in {\cal S}(x) \land rk(x,s)=0 \to y\in {\cal T}(x,s) \land N(x,s,y)=y],
\]
which is provable in $T^{1}_{2}$.

Now suppose that (\ref{eq:total1}) holds for any $b<a$ and $a>0$, and assume that 
$s\in {\cal S}(x)$ and $rk(x,s)=a$.
Let $c$ be the minimal cost of the targets in ${\cal T}(x,s)$:
\[
c:=\min\{c: \exists y<2^{d(|x|)}[ y\in {\cal T}(x,s) \land c(x,y)=c]\}
.\]
Note that ${\cal T}(x,s)$ is not empty by (\ref{eq:nPLS12}).
Pick a $y\in {\cal T}(x,s)$ such that $c(x,y)=c$.
We claim that ${\cal N}(x,s,y,y)$.
Assume ${\cal N}(x,s,y,z)$. Then by (\ref{eq:nPLS13}) with (\ref{eq:nPLS2.4}) we have 
$y,z\in {\cal T}(x,s)$, and 
either $y=z$ or $c(x,y)>c(x,z)$. The minimality of $c$ forces us to have $y=z$.

Therefore it remains to show the existence of a $u$ such that ${\cal N}(x,s,y,u)$.
By (\ref{eq:nPLS5}) 
it suffices to show the existence of a $z$ such that
${\cal N}(x,S(x,s,y),z,z)$.
By (\ref{eq:nPLS7}) we have either ${\cal N}(x,s,y,y)$ or 
$rk(x,S(x,s,y))<rk(x,s)=a$, and $S(x,s,y)\in {\cal S}(x)$ by (\ref{eq:nPLS2.1}).
Hence the Induction Hypothesis yields a $z$ such that ${\cal N}(x,S(x,s,y),z,z)$.

\hspace*{\fill} $\Box$

\section{NP search problems in $T^{2}_{2}$}\label{sect:search}

In this section we show the

\begin{theorem}\label{th:nPLST22}
If $T^{2}_{2}\vdash\forall x\exists y \, R(x,y)$ for a polynomial time computable predicate $R$, then
we can find an nPLS ${\cal P}$ in $S^{1}_{2}$ such that
\[
S^{1}_{2}\vdash y={\cal P}(x) \to R(x,U_{0}(y))
\]
for a polynomial time computable result extracting function $U_{0}$.
\end{theorem}

Before proving Theorem \ref{th:nPLST22},
we will reprove 
$\mbox{TFNP}(T^{1}_{2})\subset\mbox{PLS}$
to illuminate our proof method for searching a witness of a $\Sigma^{b}_{1}$-formula
in a derivation of bounded arithmetic.

Let $i=1,2$.
Suppose a $T^{i}_{2}$-derivation of a $\Sigma^{b}_{1}$-formula $\exists y<t(x) \, R(x,y)$ is given.
First substitute an $x$-th binary numeral $\bar{x}$ for the variable $x$, 
and unfold inference rules for induction to get another derivation essentially in the predicate logic.
Every formula occurring in the latter proof is a $\Sigma^{b}_{i}$-sentence.
The searching algorithm is so simple and canonical.
Starting with the bottom node for the end-formula $\exists y<t(\bar{x}) \, R(\bar{x},y)$
and descending the derivation tree along the Kleene-Brouwer ordering,
search a node in the derivation tree corresponding to an inference rule for introducing $\exists y<t(\bar{x}) \, R(\bar{x},y)$
and providing a true witness $\bar{n}$ such that $R(\bar{x},\bar{n})$.

Let us define derivations formally.
Extend $T^{i}_{2}$ conservatively to $T^{i}_{2}(PV)$ by adding function constants of all polynomial time computable functions.
Formulate $T^{i}_{2}(PV)$ in a one-sided sequent calculus, in which
there are extra initial sequents for axioms of (function constants for) polynomial time computable functions, and
complete induction is rendered by the inference rule with the eigenvariables $a$ and $x$
\[
\infer[(\Sigma^{b}_{i}\mbox{-ind})]
{\Gamma}
{
\Gamma, A(0)
& 
\Gamma, x\not<s_{0}, \lnot B(a,x), A(a+1)
&
\Gamma, x\not<s_{0}, \lnot B(t,x)
}
\]
for $\Sigma^{b}_{i}$-formula $A(a)\equiv \exists x<s_{0}\, B(a,x)$ where $B(a,x)$ is a $\Pi^{b}_{i-1}$-formula.
By a $\Pi^{b}_{0}$-formula or a $\Sigma^{b}_{0}$-formula we mean a literal.

By eliminating cut inferences partially, we can assume that any formula occurring in the derivation of 
$\exists y<t(x) \, R(x,y)$ is a $\Sigma^{b}_{i}$-formula.
\smallskip

Now pick a binary numeral $\bar{x}$ arbitrarily, and substitute $\bar{x}$ for the parameter $x$ in the end-formula
$\exists y<t(x) \, R(x,y)$.
Moreover unfold the inference rules ($\Sigma^{b}_{i}$\mbox{-ind}) using cut inferences:
\[
\infer[(\Sigma^{b}_{i}\mbox{-cut})]
{\Gamma}
{\cdots \Gamma,\lnot B(\bar{n}) \cdots (n<s_{0})
&
\exists x<s_{0} B(x),\Gamma
}
\]
where each $\Gamma,\lnot B(\bar{n})$ is a {\it left upper sequent\/} of the ($\Sigma^{b}_{i}$\mbox{-cut}).
$\exists x<s_{0} B(x)$ is the {\it right cut formula\/}, and each 
$\lnot B(\bar{n})$ the $n$-th {\it left cut formula\/} of the $(\Sigma^{b}_{i}\mbox{-cut})$.

In case $i=2$, unfold the combination of inference rules for bounded universal quantifiers followed by
ones for bounded existential quantifiers to the following:
\[
\infer[(\Sigma^{b}_{2})]
{\Gamma,\exists x<s_{0}\forall y<s_{1}\, L(x,y)}
{
\cdots \Gamma, \exists x<s_{0}\forall y<s_{1}\, L(x,y), L(t, \bar{n})\cdots (n<s_{1})
}
\]
where $L$ is a literal and $t$ a closed term such that $t<s_{0}$ is true.

$\Gamma,\exists x<s_{0}\forall y<s_{1}\, L(x,y)$ is the {\it lower sequent\/} of the $(\Sigma^{b}_{2})$, and
each $\Gamma, \exists x<s_{0}\forall y<s_{1}\, L(x,y), L(t, \bar{n})$ is an {\it upper sequent\/} of it.

$\exists x<s_{0}\forall y<s_{1}\, L(x,y)$ is the {\it principal formula\/} of the $(\Sigma^{b}_{2})$, and
each $L(t, \bar{n})$ is an {\it auxiliary formula\/} of it.
$t$ is the {\it witnessing term\/} of the $(\Sigma^{b}_{2})$.

This results essentially in a propositional derivation $D$ of $\exists y<t(\bar{x}) \, R(\bar{x},y)$.
Inference rules in $D$ are $(\Sigma^{b}_{i}\mbox{-cut})$, $(\Sigma^{b}_{i})$ and $(\Sigma^{b}_{1})$
\[
\infer[(\Sigma^{b}_{1})]
{\Gamma,\exists x<t\, L(x)}
{ \Gamma,\exists x<t\, L(x), L(s)\, (s<t)}
\]
where $L$ is a literal, i.e., either an equation or its negation,
 and $s$ a closed term such that $s<t$ is true.

$\Gamma,\exists x<t\, L(x)$ is the {\it lower sequent\/} of the $(\Sigma^{b}_{1})$, and
$\Gamma,\exists x<t\, L(x), L(\bar{n})$ is the {\it upper sequent\/} of the ($\Sigma^{b}_{1}$).

$\exists x<t\, L(x)$ is the {\it principal formula\/} of the $(\Sigma^{b}_{1})$, and
$L(s)$ is the {\it auxiliary formula\/} of it.
$s$ is the {\it witnessing term\/} of the $(\Sigma^{b}_{1})$.

Note the fact:
\begin{equation}\label{eq:upperlower}
\mbox{Each upper sequent contains its lower sequent.}
\end{equation}

There occurs no free variable in $D$.
Initial sequents in $D$ are
\[
\Gamma, L
\]
for true literals $L$.

Note that we can calculate the values of closed terms in $D$, in polynomial time
since we can assume that we concern only a finite number of polynomial time computable functions, and
each closed term in $D$ has a constant depth.

It is easy to see that  for a polynomial $d$,
the size of $D$(the number of occurrences of symbols in $D$) is bounded by $2^{d(|x|)}$, and
the depth of $D$ is bounded by $d(|x|)$.

Observe that each sequent in $D$ consists of $\Sigma^{b}_{i}$, $\Sigma^{b}_{1}$ and $\Sigma^{b}_{0}$(literals) 
sentences:
\[
\Gamma=\Gamma(\Sigma^{b}_{i}),\Gamma(\Sigma^{b}_{1}),\Gamma(\Sigma^{b}_{0})
\]
where $\Gamma(\Sigma^{b}_{k})\subseteq \Sigma^{b}_{k}$ for $k=1,2$.

To extract informations in $D$, let us specify what is a derivation of this 'propositional' calculus.
Let $T(D)$ denote a naked tree of $D$.
Each node $\sigma\in T(D)$ is a finite sequence of natural numbers.
The empty sequence $\emptyset$ is the root.
$Seq_{\sigma}$ for $\sigma\in D(T)$ denotes the sequent situated at the node $\sigma$.
It is assumed that we can calculate, in polynomial time from $\sigma$, the principal formula, auxiliary formulas, witnessing term 
when $Seq_{\sigma}$ is a lower sequent of a $(\Sigma^{b}_{k})\, (k=1,2)$,
and the cut formula when $Seq_{\sigma}$ is an upper sequent of a $(\Sigma^{b}_{i}\mbox{-cut})$.

Let $<_{KB}$ denote the Kleene-Brouwer ordering on the tree $T(D)$.
It is easy to see the existence of a polynomial time computable function 
$kb: D\ni \sigma\mapsto kb(x,\sigma)$, which is
compatible with the Kleene-Brouwer ordering $<_{KB}$ on $T(D)$:
\begin{equation}\label{eq:KB}
\sigma,\tau\in T(D) \,\&\, \sigma<_{KB}\tau \Rightarrow kb(x,\sigma)<kb(x,\tau)
\end{equation}

To arithmetize formal objects such as terms, formulas, sequents and derivation,
 we assume a feasible encoding of finite sequences of natural numbers.
${}^{<\omega}\omega$ denotes the set of sequence numbers, $Len(\sigma)$ the length of the sequence coded by the number 
$\sigma\in {}^{<\omega}\omega$, $(\sigma)_{k}$ the $k$-th entry in $\sigma$ for $k<Len(\sigma)$, $\sigma*\tau$ concatenated sequence, and $\sigma\subseteq\tau$ means that $\sigma$ is an initial segent of $\tau$, etc.
All of these are $\Delta^{b}_{1}$-definable in $S^{1}_{2}$, and the bounded arithmetic $S^{1}_{2}$
proves the elementary facts on them, cf. \cite{Buss}.

\subsection{$\mbox{TFNP}(T^{1}_{2})\subset\mbox{PLS}$}

In this subsection consider the case $i=1$, and we define a PLS problem from the derivation $D$ of $\exists y<t(\bar{x}) \, R(\bar{x},y)$.
Inference rules in $D$ are $(\Sigma^{b}_{1}\mbox{-cut})$ and $(\Sigma^{b}_{1})$.
The algorithm is based on the following simple observation.

Let $\sigma$ be a node in the tree $T(D)$ such that
 \begin{equation}\label{eq:targetcnd}
 Seq_{\sigma} \mbox{ contains no true literal.}
 \end{equation}
Let us view such a node $\sigma$ as a search problem in searching a witness for
a $\Sigma^{b}_{1}$-formula in $Seq_{\sigma}(\Sigma^{b}_{1})$.

\begin{proposition}\label{prp:sigma1}
Let $\sigma$ be a node in the tree $T(D)$ enjoying the condition (\ref{eq:targetcnd}).
Then there exists a $\tau\supseteq\sigma$ on the rightmost branch in the upper part of $\sigma$ such that
$Seq_{\tau}$ is a lower sequent of a $(\Sigma^{b}_{1})$ whose auxiliary formula is true.
Moreover the lowest such node $\tau=t(x,\sigma)$ is polynomial time computable.
\end{proposition}
{\it Proof}.\hspace{2mm}
Consider the rightmost branch in the upper part of $\sigma$.
Its top contains a true literal $L$, and the literal disappears before reaching to $\sigma$
 by the condition (\ref{eq:targetcnd}).
The vanishing point has to be a $(\Sigma^{b}_{1})$ with its auxiliary formula $L(s)\equiv L$:
\[
\infer[(\Sigma^{b}_{1})]
{\tau: \Gamma,\exists x<t\, L(x)}
{ \Gamma,\exists x<t\, L(x), L(s)}
\]
since the rightmost branch does not pass through any left upper sequent of a 
$(\Sigma^{b}_{1}\mbox{-cut})$.
This shows the existence of a node $\tau$.
The lowest such node $\tau=t(x,\sigma)$ is polynomial time computable
since the depth of the tree $T(D)$ is bounded by a polynomial $d(|x|)$.
\hspace*{\fill} $\Box$

There are three cases to consider according to the vanishing point $\lambda$ 
of the principal formula $A\equiv\exists x<t\, L(x)$ of the lowest 
$\tau: (\Sigma^{b}_{1})\, (\tau=t(x,\sigma))$.
Let $n$ denote the value of its witnessing term $s$.

\begin{enumerate}
\item
$\lambda=\emptyset$:
 Namely $A$ is the end-formula $\exists y<t(\bar{x}) \, R(\bar{x},y)$.
Then we are done, and the $n$ is a witness sought.

Otherwise $A$ is a right cut formula of a $\lambda:(\Sigma^{b}_{1}\mbox{-cut})$.

\item
$\lambda\supseteq\sigma$:
Then $A\not\in \Gamma_{1}\equiv Seq_{\sigma}$:
 \[
 \infer*
 {\sigma: \Gamma_{1}}
{
 \infer[(\Sigma^{b}_{1}\mbox{-cut})]
 {\lambda: \Gamma}
  {
   \cdots \kappa : \Gamma,\lnot L(\bar{n}) \cdots 
  &
   \infer*{A,\Gamma}
   {
    \infer[(\Sigma^{b}_{1})]
    {\tau : A, \Delta_{0}}{ A,\Delta_{0}, L(\bar{n})}
   }
  }
 }
\]

\item
$\lambda\subset\sigma$:
Then $A\in\Gamma_{1}$:

 \[
 \infer[(\Sigma^{b}_{1}\mbox{-cut})]
  {\lambda: \Gamma}
{
  \cdots \kappa : \Gamma,\lnot L(\bar{n}) \cdots
  &
  \infer*{A,\Gamma}
   {
    \infer*{ \sigma: A, \Gamma_{1}}
     {
      \infer[(\Sigma^{b}_{1})]
      {\tau : A, \Delta_{0}}{A,\Delta_{0}, L(\bar{n})}
      }
  }
}
\]
\item
In the latter two cases, $\kappa$ enjoys the condition (\ref{eq:targetcnd}).
Since $\lnot L(\bar{n})$ is a false literal, it suffices to see that $\Gamma=Seq_{\lambda}$ contains no true literal in the above figures.

First let $\lambda\supseteq\sigma$, and consider $\rho$ with $\lambda\supseteq\rho\supset\sigma$.
Then $Seq_{\rho}$ is either a right upper sequent of a $(\Sigma^{b}_{1}\mbox{-cut})$
or an upper sequent of a $(\Sigma^{b}_{1})$ with a false auxiliary formula.
Therefore if its lower sequent enjoys (\ref{eq:targetcnd}), then so does the upper $\rho$.

Next assume $\lambda\subset\sigma$.
Then by (\ref{eq:upperlower})
 $Seq_{\lambda}=\Gamma\subseteq\Gamma_{1}=Seq_{\sigma}$.
\end{enumerate}

Therefore $\kappa$ is another searching problem.
Here is another simple observation.

\begin{proposition}\label{prp:KB}
In the above figures, 
\[
\kappa<_{KB}\sigma
\]
and $\sigma\mapsto\kappa=f(x,\sigma)$ is a polynomial time computable map.
\end{proposition}

Now let us define a PLS for searching a witness of the end-formula $\exists y<t(\bar{x}) \, R(\bar{x},y)$.

\begin{enumerate}
\item
The set $F(x)$ of feasible points consists in $\sigma\in T(D)$ such that
either $\sigma=\emptyset$ or $Seq_{\sigma}$ is a left upper sequent of a $(\Sigma^{b}_{1}\mbox{-cut})$ enjoying the condition (\ref{eq:targetcnd}).
\item
The initial point $i(x)=\emptyset$.
\item
The neighborhood function $N(x,\sigma)=\kappa$ for $\sigma\in F(x)$ is defined
by $N(x,\sigma)=f(x,\sigma)$
for the function $f$ in Proposition \ref{prp:KB} except the first case $\lambda=\emptyset$.
This means that $\tau=t(x,\sigma)$ in Proposition \ref{prp:sigma1} is a lower sequent
of a $(\Sigma^{b}_{1})$ whose auxiliary formula is a true literal $R(\bar{x},\bar{n})$.
Then $N(x,\sigma)=\sigma$.
\item
The cost function $c(x,\sigma)=kb(x,\sigma)$ in (\ref{eq:KB}).
\end{enumerate}
Then we can extract a witness $n$ such that $R(\bar{x},\bar{n})$ is true
from any solution $\sigma$ of the PLS.

\subsection{Proof of Theorem \ref{th:nPLST22}}
Now let $D$ be a derivation of $\Sigma^{b}_{1}$-sentence $\exists y<t(\bar{x}) R(\bar{x},y)$,
which arises from a derivation in $T^{2}_{2}$, and define an nPLS problem 
to extract a witness.
Inference rules in $D$ are $(\Sigma^{b}_{2}\mbox{-cut})$, $(\Sigma^{b}_{1})$ and $(\Sigma^{b}_{2})$.
Again a node $\sigma$ enjoying the condition (\ref{eq:targetcnd}) denotes
 a search problem in searching a witness for
a $\Sigma^{b}_{1}$-formula in $Seq_{\sigma}(\Sigma^{b}_{1})$.

Proposition \ref{prp:sigma1} is modified as follows.
\begin{proposition}\label{prp:sigma2}
Let $\sigma$ be a node in the tree $T(D)$ enjoying the condition (\ref{eq:targetcnd}).
Then there exists a $\tau\supseteq\sigma$ on the rightmost branch in the upper part of $\sigma$ such that
either $Seq_{\tau}$ is a lower sequent of a $(\Sigma^{b}_{1})$ whose auxiliary formula is true,
or $Seq_{\tau}$ is a lower sequent of a $(\Sigma^{b}_{2})$.
Moreover the lowest such node $\tau=t(x,\sigma)$ is polynomial time computable.
\end{proposition}

There are three cases to consider according to the principal formula $A\equiv\exists x<t\, B(x)$ of the lowest 
$\tau=t(x,\sigma)$.
Let $n$ denote the value of its witnessing term $s$, and $\lambda$ the vanishing point of the principal formula.

\begin{enumerate}
\item
$A$ is a $\Sigma^{b}_{2}$-sentence:
Then $\lambda$ is a $(\Sigma^{b}_{2}\mbox{-cut})$, and $A$ is its right cut formula.
Let $B(x)\equiv\forall y<s L(x,y)$ with a literal $L$.
 \[
 \infer[(\Sigma^{b}_{2}\mbox{-cut})]
 {\lambda: \Gamma}
  {
   \cdots \kappa : \Gamma,\lnot B(\bar{n}) \cdots 
  &
   \infer*{A,\Gamma}
   {
    \infer[(\Sigma^{b}_{2})]
    {\tau : A, \Delta_{0}}{ \cdots A,\Delta_{0}, L(\bar{n},\bar{m}) \cdots}
   }
  }
\]

We claim that $\kappa$ enjoys the condition (\ref{eq:targetcnd}).
Since $\lnot B(\bar{n})$ is a $\Sigma^{b}_{1}$-sentence, it suffices to see that $\Gamma=Seq_{\lambda}$ contains no true literal in the figure.

By (\ref{eq:upperlower}) we can assume $\lambda\supseteq\sigma$.
Let $\rho$ with $\lambda\supseteq\rho\supset\sigma$.
Then $Seq_{\rho}$ is either a right upper sequent of a $(\Sigma^{b}_{2}\mbox{-cut})$
or an upper sequent of a $(\Sigma^{b}_{1})$ with a false auxiliary formula.
Hence if its lower sequent enjoys (\ref{eq:targetcnd}), then so does the upper $\rho$.

Therefore $\kappa$ is another searching problem, cf. the definition of $S(x,\sigma,\tau)$ below.
In the above figure, we have
\[
\kappa<_{KB}\sigma
\]
and $\sigma\mapsto\kappa=f(x,\sigma)$ is a polynomial time computable map, cf. Proposition \ref{prp:KB}.

\item
 $A$ is the end-formula $\exists y<t(\bar{x}) \, R(\bar{x},y)$.
Then we are done, and the $n$ is a witness sought.

\item
Otherwise:
 $A$ is a left cut formula $\lnot B(\bar{p})$ of a $\lambda:(\Sigma^{b}_{2}\mbox{-cut})$.
Let $A\equiv\exists y<s \lnot L(\bar{p},y)$ with a literal $L$.
We have $\lambda*\langle p\rangle\subseteq\sigma$.

Consider the case when $\lambda*\langle p\rangle=f(x,\sigma_{0})$ for a problem $\sigma_{0}$ with
$\lambda*\langle p\rangle<_{KB}\sigma_{0}$:
 \[
    \infer[(\Sigma^{b}_{2}\mbox{-cut})]
    {\lambda: \Gamma}
    {
    \infer*{\cdots \lambda*\langle p\rangle:\Gamma,\lnot B(\bar{p}) \cdots(p<t_{0})}
     {
      \infer[(\Sigma^{b}_{1})]
      {\tau : A, \Delta_{0}}{A,\Delta_{0}, \lnot L(\bar{p},\bar{n})}
      }
    &
       \infer*{\exists x<t_{0} B(x),\Gamma}
       {
        \infer[(\Sigma^{b}_{2})]{\tau_{0}: \exists x<t_{0} B(x),\Gamma_{0}}
        {\cdots \tau_{0}*\langle n\rangle: \exists x<t_{0} B(x),\Gamma_{0}, L(\bar{p},\bar{n})\cdots(n<s)}
       }
    }
\]
where $\tau_{0}=t(x,\sigma_{0})$ and $L(\bar{p},\bar{n})$ is a false literal.

Then the solution $n$ of the problem $\lambda*\langle p\rangle:\Gamma,\lnot B(\bar{p})$
tells the problem $\sigma_{0}$ where to proceed.
Namely go to the $n$-th path $\tau_{0}*\langle n\rangle: \exists x<t_{0} B(x),\Gamma_{0}, L(\bar{p},\bar{n})$,
and then climb up to the node $t(x,\tau_{0}*\langle n\rangle)$, cf. the definition of $U(x,\sigma,\tau, \rho)$ below.

\end{enumerate}

Now let us define an nPLS problem 
\[
{\cal P}=\langle d; {\cal S},{\cal T},{\cal N};N,  i,t,c,S,U,rk\rangle
.\]

\begin{enumerate}
\item
$\sigma\in{\cal S}(x)$ iff $\sigma\in T(D)$ and the sequent $Seq_{\sigma}$ is one of the following:
 \begin{enumerate}
 \item
 $Seq_{\sigma}$ is the endsequent, i.e., $\sigma=\emptyset$.
 \item
 $Seq_{\sigma}$ is one of a left upper sequent of a $(\Sigma^{b}_{2}\mbox{-cut})$ such that
\begin{eqnarray}
 &&
 \forall \tau[\emptyset\subseteq\tau\subset\sigma \to Seq_{\tau} 
 \mbox{ is not a lower sequent of any } (\Sigma^{b}_{1}) 
  \nonumber
  \\
 &&
 \mbox{ whose auxiliary formula is a true literal}]
 \label{eq:sourcecnd}
\end{eqnarray}
 \end{enumerate}
Then $\sigma$ enjoys the condition (\ref{eq:targetcnd}).
 
Each $\sigma\in{\cal S}(x)$ corresponds to a search problem in searching a witness for
a $\Sigma^{b}_{1}$-formula in $Seq_{\sigma}(\Sigma^{b}_{1})$.

\item
$\tau\in{\cal T}(x,\sigma)$ iff $\sigma,\tau\in T(D)$, $\sigma\in{\cal S}(x)$, and
the sequent $Seq_{\tau}$ is one of the following two: 
  \begin{enumerate}
 \item
 $\sigma\subseteq\tau$, $Seq_{\tau}$ is a lower sequent of a $(\Sigma^{b}_{2})$, and
 \begin{equation}\label{eq:target2}
 \mbox{there is no left upper sequent of any } (\Sigma^{b}_{2}\mbox{-cut}) \mbox{ between } \sigma
 \mbox{ and } \tau
 \end{equation}
This means that the path from $\sigma$ to $\tau$ is in the rightmost branch of the upper part of $\sigma$.

  \item
 $Seq_{\tau}$ is a lower sequent of a $(\Sigma^{b}_{1})$ such that its auxiliary formula is true, and
 its principal formula is in $Seq_{\sigma}$.
 \end{enumerate}

Moreover in each case $Seq_{\tau}$ has to enjoy the condition (\ref{eq:targetcnd}).

\item
For $\tau\in {\cal T}(x,\sigma)$, $c(x,\tau)$ denotes the depth of $\tau$ in the finite tree $T(D)$
if $Seq_{\tau}$ is a lower sequent of a $(\Sigma^{b}_{2})$.
Otherwise put $c(x,s)=0$.

\item
$rk(x,\sigma)=kb(x,\sigma)$ in (\ref{eq:KB}).

\item
$i(x)=\emptyset$.

\item
For $\sigma\in {\cal S}(x)$, $t(x,\sigma)$ denotes the node $\tau\in{\cal T}(x,\sigma)$ which is
 the lower sequent of the lowest $(\Sigma^{b}_{i})\, (i=1,2)$
 in the rightmost branch of the upper part of $\sigma$, cf. Proposition \ref{prp:sigma2}.

\item

For $\tau,\rho\in {\cal T}(x,\sigma)$, 
${\cal N}(x,\sigma,\tau,\rho)$ iff  
one of the following two holds: 
  \begin{enumerate}
 \item
$Seq_{\tau}$ is a lower sequent of a $(\Sigma^{b}_{2})$, and
if $Seq_{\rho}$ is a lower sequent of a $(\Sigma^{b}_{2})$, then $\tau\subset\rho$.

  \item
 $Seq_{\tau}$ is a lower sequent of a $(\Sigma^{b}_{1})$, and $\tau=\rho$.
 \end{enumerate}

Note that if ${\cal N}(x,\sigma,\tau,\tau)$, then the witnessing term of the $(\Sigma^{b}_{1})$
whose lower sequent is $Seq_{\tau}$, yields a witness for the search problem $\sigma$.

\item
For $rk(x,\sigma)=0$ and $\tau\in{\cal T}(x,\sigma)$, $N(x,\sigma,\tau)=\rho$ is defined as follows.
First if $Seq_{\tau}$ is a lower sequent of a $(\Sigma^{b}_{1})$, then $\rho:=\tau$.
Next suppose that $Seq_{\tau}$ is a lower sequent of a $(\Sigma^{b}_{2})$ with its principal formula $A$.
Consider the $(\Sigma^{b}_{2}\mbox{-cut})$ below $\tau$ whose right cut formula is $A$, and let $\kappa$ denote
a left upper sequent of the $(\Sigma^{b}_{2}\mbox{-cut})$.
Then $\kappa\not\in{\cal S}(x)$ since $\kappa<_{KB}\sigma$ and $rk(x,\sigma)=0$.
This means by (\ref{eq:sourcecnd}) that there is a $\tau$ such that
$\emptyset\subseteq\tau\subset\kappa$ and
$Seq_{\tau}$ is a lower sequent of a $(\Sigma^{b}_{1})$ whose auxiliary formula is a true literal.
This is not the case since (\ref{eq:sourcecnd}) for $\sigma\in{\cal S}(x)$,
(\ref{eq:targetcnd}) for $\tau\in{\cal T}(x,\sigma)$,
and $\lambda$ is comparable with $\sigma$ with respect to
 the tree ordering $\subseteq$ for $\kappa=\lambda*\langle n\rangle$.

\item
For $\tau\in {\cal T}(x,\sigma)$,
 $S(x,\sigma,\tau)\in {\cal S}(x)$ is defined as follows.

 Suppose $\lnot{\cal N}(x,\sigma,\tau,\tau)$.
 This means that $\tau$ is a lower sequent of a $(\Sigma^{b}_{2})$.
 Define $S(x,\sigma,\tau)=f(x,\sigma)$ in Proposition \ref{prp:KB}.
 
Let $A\equiv \exists x<s_{0}\forall y<s_{1}\, L(x,y)$
  be the principal formula of the $(\Sigma^{b}_{2})$,
 $t$ its witnessing term with the value $n$.
 
  Note that $S(x,\sigma,\tau)<_{KB}\sigma$ for (\ref{eq:nPLS7}) since
 $\sigma\subseteq\tau$ and the vanishing point $\lambda$ of $A$  is comparable with $\sigma$ with respect to $\subseteq$.
 
 Let $\kappa=S(x,\sigma,\tau)$.
 First the case when $\sigma\subseteq\lambda$:
 \[
 \infer*
 {\sigma: \Gamma_{1}}
{
 \infer[(\Sigma^{b}_{2}\mbox{-cut})]
 {\lambda: \Gamma}
  {
   \cdots \kappa : \Gamma,\lnot\forall y<s_{1}\, L(\bar{n},y) \cdots 
  &
   \infer*{A,\Gamma}
   {
    \infer[(\Sigma^{b}_{2})]
    {\tau : A, \Delta_{0}}{ \cdots A,\Delta_{0}, L(\bar{n},\bar{k})\cdots\, (k<s_{1})}
   }
  }
}
\]
Note that the condition (\ref{eq:sourcecnd}) for $\kappa\in{\cal S}(x)$ is enjoyed 
by (\ref{eq:targetcnd}) for $\tau\in{\cal T}(x,\sigma)$ and (\ref{eq:upperlower}).

Second the case when $\lambda\subset\sigma$:

 \[
 \infer[(\Sigma^{b}_{2}\mbox{-cut})]
  {\lambda: \Gamma}
  {
  \cdots \kappa : \Gamma,\lnot\forall y<s_{1}\, L(\bar{n},y) \cdots
  &
  \hspace{-20mm}
  \infer*{A,\Gamma}
    {
    \infer[(\Sigma^{b}_{2}\mbox{-cut})]
    {A,\Gamma_{1}}
    {
    \infer*{\cdots \sigma: A, \Gamma_{1},\lnot B(\bar{p}) \cdots}
     {
      \infer[(\Sigma^{b}_{2})]
      {\tau : A, \Delta_{0}}{ \cdots A,\Delta_{0}, L(\bar{n},\bar{k})\cdots\, (k<s_{1})}
      }
    &
    \hspace{-7mm}
    \exists x<t_{0} B(x),\Gamma_{1}
    }
  }
}
\]

\item

For 
${\cal N}(x,S(x,\sigma,\tau),\rho,\rho)$,
$U(x,\sigma,\tau, \rho)\in {\cal T}(x,\sigma)$ is defined as follows.

 Let $n$ denote the value of the witnessing term of the $(\Sigma^{b}_{1})$ 
 whose lower sequent is $\rho$.
$\tau$ is the lower sequent of a $(\Sigma^{b}_{2})$:
\[
\infer[(\Sigma^{b}_{2})]
{\tau: \Gamma,\exists x<s_{0}\forall y<s_{1}\, L(x,y)}
{
\cdots \tau*\langle n\rangle : \Gamma, \exists x<s_{0}\forall y<s_{1}\, L(x,y), L(t, \bar{n})\cdots (n<s_{1})
}
\]
Let $\kappa=S(x,\sigma,\tau)$.
There are two cases to consider.
 \begin{enumerate}
  \item
 The case when the principal formula of $\rho$
  is in $Seq_{\kappa}$, but not in the lower sequent of the $(\Sigma^{b}_{2}\mbox{-cut})$:

Let $U(x,\sigma,\tau, \rho)$ denote the rightmost and lowest node $\xi\in {\cal T}(x,\sigma)$ 
such that $\tau*\langle n\rangle \subseteq \xi=t(x,\tau*\langle n\rangle)$.
Such a node exists in the rightmost branch of the upper part of $\tau*\langle n\rangle$
since the literal $L(t,\bar{n})$ is false, cf. Proposition \ref{prp:sigma2}.

 \item
  Otherwise:

 From (\ref{eq:target2}) (when $\sigma\subset \kappa$) and the fact (\ref{eq:upperlower})  
 (when $\lambda\subset \sigma$ for $\lambda*\langle n\rangle =\kappa$),
 we see that $\rho$ provides a solution for a search problem for 
$\sigma$.
 
 Put $U(x,\sigma,\tau, \rho)=\rho$.

 \end{enumerate}

We see that (\ref{eq:nPLS5})  for ${\cal N}(x,\sigma,\tau,U(x,\sigma,\tau, \rho))$ is enjoyed in each case.

\end{enumerate}

Assume ${\cal N}(x,i(x),\tau,\tau)$.
Then $\tau$ is a $(\Sigma^{b}_{1})$ witnessing the end-formula of $D$.
Therefore the value $U_{0}(\tau)$ of its witnessing term realizes Theorem \ref{th:nPLST22}.
\\ \smallskip \\

For future works, the approach in this note could be extended to characterize
 the $\Sigma^{b}_{1}$-definable functions in $T^{k}_{2}$ for any $k\geq 2$
 by introducing classes of search problems of higher order PLS.




\end{document}